\documentclass[12pt]{article}
\usepackage{amsmath}
\usepackage{graphicx}

\usepackage{tikz}
\usepackage{amssymb} 
\usepackage{colortbl}
\newtheorem{thm}{Theorem}
\newtheorem{lemma}[thm]{Lemma}

\newtheorem{prop}[thm]{Proposition}

\begin{document}
\title{The  Infinity-Laplacian in Smooth Convex Domains and in a Square}
\author{Karl K. Brustad, Erik Lindgren, Peter Lindqvist} 
\date{\today}

\maketitle
\qquad \qquad  \emph{Dedicated to Giuseppe Mingione on his fiftieth  birthday}

\medskip
{\small \textsc{Abstract:} \textsf{We extend some theorems for the Infinity-Ground State and for the Infinity-Potential, known for convex polygons, to other domains  in the plane, by applying Alexandroff's method to the curved boundary. A recent \emph{explicit} solution disproves a conjecture.}

\bigskip 
\noindent {\small \textsf{AMS Classification 2000}: 35J65, 35J94, 35P30, 49N60.}


\section{Introduction}

The $\infty$\,-\,Laplace Operator
$$\Delta_{\infty}u\,\equiv\,\sum_{i,j}\frac{\partial u}{\partial x_{i}}\, \frac{\partial u}{\partial x_{j}}\,\frac{\partial^2 u}{\partial x_i\partial x_j}$$
is formally the limit of the $p$\,-\,Laplace Operator
$$\Delta_p u\,\equiv\,\nabla\!\cdot\!(|\nabla u|^{p-2}\nabla u)\quad\text{as}\quad p \to \infty.$$
The two-dimensional equation
$$\Bigl(\frac{\partial u}{\partial x_1}\Bigr)^{\!2}\frac{\partial ^{\,2} u}{\partial x_1^{\,2}} +2\,\frac{\partial u}{\partial x_1} \frac{\partial u}{\partial x_2}\frac{\partial^{\,2}u}{\partial x_1 \partial x_2} + \Bigl(\frac{\partial u}{\partial x_2}\Bigr)^{\!2}\frac{\partial {\,^2} u}{\partial x_2^{\,2}}\,\,=\,\,0$$
was introduced by G. Aronsson in 1967 as a tool to provide optimal\mbox{ Lipschitz} extensions, cf. [A1], [A2]. It has been intensively studied ever since; some highlights are
\begin{itemize}
\item Viscosity solutions for $\Delta_{\infty}$ were introduced by T. Bhattacharya, E. DiBenedetto, and J. Manfredi in [BDM].
\item R. Jensen proved uniqueness in [J].
\item Differentiability was proved by O. Savin and L. Evans. See [S1] and [ES].
\item The connexion with stochastic game theory (``Tug-of-War'') was discovered by Y. Peres, O. Schramm, S. Sheffield, and D. Wilson, cf. [PSSW].
\end{itemize}

By examples we shall shed some light on two  problems in the plane. The problems  are related but not identical. In convex domains so  similar methods often work for both problems that it is optimal to treat them simultaneously. ---We remark that the main difficulty is the lack of second derivatives. The  solutions are to be interpreted as $\emph{viscosity solutions}$. The reader may consult [Ko] and [CIL]. 

The first one is the boundary value problem
\begin{equation}\label{ring}
  \begin{cases}
    \Delta_{\infty}u&=\,0\quad \text{in} \quad G\\
   \quad u&=\,0 \quad \text{on} \quad \partial \Omega\\
  \quad  u&=\,1 \quad \text{on} \quad \partial K
  \end{cases}
  \end{equation}
    in \emph{a convex ring} $G = \Omega \setminus K$, where $\Omega$ is a bounded convex domain in $\mathbb{R}^2$ and $K \subset \Omega $ is a closed convex set. The unique solution, say $u_{\infty}$, belongs to $C(\overline{G})$ and always takes the boundary values; $K$ is often only an isolated point. We say that $u_{\infty}$ is the $\infty$-potential. In [L] the term ``capacitary function'' is used. 

    The second object is the $\infty$\,-\,Eigenvalue Problem
    \begin{equation}\label{Lambda}\begin{cases}
        &\mathrm{max}\left\{\Lambda-\dfrac{|\nabla v|}{v},\,\Delta_{\infty}v \right\}\,=\,0\quad \text{in}\quad \Omega\\
       & v \in C(\overline{\Omega}),\,\, v|_{\partial \Omega}\,=\,0,\,\, v > 0.
      \end{cases}
    \end{equation}
    Solutions are called $\infty$\,-\,Ground States. Problem (\ref{Lambda}) is the asymptotic limit as $p\to\infty$ of the equation
    \begin{equation}\label{450}
      \nabla\!\cdot\!\bigl(|\nabla v_p|^{p-2}\nabla v_p\bigr)\,+\, \lambda_p|v_p|^{p-2}v_p\,=\,0, 
      \end{equation}
    where $v_p \in W^{1,p}_0(\Omega),\,\, v_p > 0.$ Problem (\ref{Lambda}) has a solution if and only if
    $$\Lambda\,=\, \lim_{p\to \infty} \sqrt[\leftroot{-2}\uproot{3}p]{\lambda_p},$$
    but uniqueness (= simplicity of $\Lambda$) is not known to hold even in convex domains.\footnote{In more general domains this fails, cf. [HSY].}

    We shall restrict ourselves to those $\infty$\,-\,Ground States that come as limits of sequences of solutions to (\ref{450}). Such a limit $v$ has the advantage that $\log v$ is concave. This valuable property is the reason for why we prefer these so-called \emph{variational} $\infty$\,-\,Ground States. See [JLM1] and [Y]. 

The main achievement of this paper is to complement our study in [LL2] and [LL3]. There the crucial assumption that $|\nabla u_{\infty}| $  (or $|\nabla v_{\infty}|$) has only a finite number of maxima and minima on the boundary $\partial \Omega$ was properly verified merely for convex polygons. Smooth domains were out of reach.\footnote{We think that every bounded convex domain with $C^3$ - boundary will do.} Our contribution now is to provide a class of explicit smooth domains having the desired property: for example, the ellipse is included. We are grateful to B. Kawohl, who informed us about [Ka] and suggested Alexandroff's Moving Plane Method. 
    
Despite sharing similar properties, solutions of  \eqref{ring} and \eqref{Lambda} may not coincide even under the (necessary) condition that $K$ be chosen as the \emph{High Ridge}
\begin{equation}\label{High}
  K\,=\,\{x \in \Omega|\, \mathrm{dist}(x,\partial \Omega)\,=\,R\},\quad
  R\,=\,\max_{x\in \Omega}\,\mathrm{dist}(x,\partial\Omega).
\end{equation}
(Here $R$
is the radius of the largest inscribed ball in $\Omega$.)
Nevertheless, it is shown in Theorem 3.3 in [Y], that in a certain class of domains, which  includes the  stadium--like domains, the distance function is the unique solution to both \eqref{ring} and \eqref{Lambda}. In general, coincidence is a difficult problem.

We conclude the work by noticing a recent result in a punctured square: the ring  domain is \emph{a square with its center removed}. It has been predicted that the $\infty$-Potential would coincide with the $\infty$\,-\,Ground State, cf. [JLM2]. Brustad's explicit formula in [B] reveals that \emph{the functions do not coincide}. However, the maximal difference between the functions is $< 10^{-3}$ for a square of area $4$, according to numerical calculations in [BBT].

     \paragraph{Acknowledgements:}  We thank Bernd Kawohl for his valuable piece of advice. 
     This work was done while the authors during the fall of 2022  were participating in the research program ``Geometric Aspects of Nonlinear Partial Differential Equations'' at Institut Mittag-Leffler. It was  supported by the Swedish Research Council under grant no. 2016-06596,
    E. L. is supported by the Swedish Research Council, grant no. 2017-03736. 

    \section{Preliminaries}\label{sec:prel}

    We use standard notation. Here $\Omega \subset \mathbb{R}^2$ will always denote a bounded convex domain with smooth boundary, say at least of class $C^{2,\alpha}$.  We shall denote the $\infty$\,-\,Ground State by $v_{\infty}$, and the solution of (\ref{450}) by $v_p$. Analogously, the solution of (\ref{ring}) is $u_{\infty}$, and $u_p$ will denote the solution of problem (\ref{pharmappr}), see Section \ref{sec:ppot}. Thus
    $$v_{\infty}\,=\,\lim_{p\to \infty} v_p,\qquad u_{\infty}\,=\,\lim_{p\to \infty} u_p,$$ perhaps via a subsequence.  We will also use the normalization 
    \begin{equation}\label{eq:normal}
      \underset{x \in \Omega}{\max}\,v_\infty\,=\,1.
    \end{equation}
    According to Theorem 2.4 in [Yu], the High Ridge (see \eqref{High}) is also the set where $v_\infty$ attains its maximum $1$.
    
    \paragraph{The contact set} The contact set  $\Upsilon$ plays a central role. To define it, following Y. Yu  we use the operator
    $$ S^-(x)\,=\,\lim_{r \to 0}\,\Bigl \{- \min_{y \in \partial B(x,r)}\frac{v_{\infty}(y)-v_{\infty}(x)}{r}\Bigr \}.$$
    According to Theorem 3.6 in [Yu], $S^-$ is continuous in $\Omega$ and at points of differentiability  $S^-(x) = |\nabla v_{\infty}(x)|.$ The \emph{contact set} $$ \Upsilon\,=\,\bigl\{x\in \Omega\,| \, S^-(x)\,=\,\Lambda v_\infty(x)\bigr\}$$
is closed and has zero area, see Corollary 3.7 in [Yu]. In addition, by Corollary 7 in [LL3], $\Upsilon$ does not reach $\partial \Omega$. The set where $v_{\infty}$ is not differentiable is contained in $\Upsilon$, cf. Lemma 3.5 in [Yu]. In the open set $\Omega\setminus \Upsilon$, $v_\infty$ is $C^1$ and $\Delta_\infty v_\infty=0$, by Theorem 3.1 in [Yu]. 

   \paragraph{Streamlines}  For the benefit of the reader, we describe the role of the streamlines. All this can be found in [LL2] and [LL3]. We begin with the $\infty$-potential $u_{\infty}$. From each boundary point $\xi \in \partial \Omega$ a unique streamline $\boldsymbol\alpha = \boldsymbol\alpha(t)$ starts and reaches $K$ in finite time:
    $$\frac{d\boldsymbol\alpha(t)}{dt}\,=\,+\,\nabla u_{\infty}(\boldsymbol\alpha(t)),\qquad \boldsymbol\alpha(0)\,=\,\xi\,;\quad 0\leq t\leq T.$$
    It may meet and join other streamlines, but streamlines do not cross. Suppose now that\footnote{We do not know of any convex plane domain for which this is not valid!}
    \begin{quote} {\sf Along $\partial \Omega$ the speed $ |\nabla u_{\infty}(\xi)|$ has only a finite number of local minima.}
    \end{quote}
    More exactly we allow a finite number of strict local minimum points and a finite number of boundary arcs along which strict minima are obtained.
    It is problematic to deduce this from the shape of the domain. Polygons were treated in [LL2], and for a family of smooth domains the above assumption will be verified in Section 5. --- This is our achievement in the present work.

    The streamline starting at a  \emph{strict} local minimum point is called \emph{an attracting streamline}. If a local minimum is attained along a whole closed boundary arc (there $|\nabla u_{\infty}|$ is constant), then this minimum produces \emph{two} attracting streamlines: the two streamlines emerging from the endpoints of the arc. The attracting streamlines are special, indeed. We cite the main theorem from [LL2].

    \begin{thm}\label{LLU} Let $\boldsymbol\alpha$ be a streamline of $u_{\infty}$ that is not an attracting one. Then it cannot meet any other streamline before it either meets (and joins) an attracting streamline or reaches $K$. The speed $|\nabla u_{\infty}(\boldsymbol\alpha(t))|$ is constant along $\boldsymbol\alpha$ until it joins an attracting streamline, after which the speed is non-decreasing.
    \end{thm}

    The corresponding theorem for the streamlines of $v_{\infty}$ is similar; we only  have to replace $K$ by the \emph{High Ridge} of $\Omega$; see (\ref{High}).
       The High Ridge is also the set of points at which  $v_{\infty}$ attains its maximum.
    Suppose that $\boldsymbol \gamma_1,\boldsymbol \gamma_2,\cdots,\boldsymbol \gamma_N$ are the attracting streamlines. It is spectacular that outside the closed set
    $$\Gamma\,=\, \boldsymbol \gamma_1 \cup\boldsymbol \gamma_2 \cup \cdots \cup \boldsymbol \gamma_N   $$
    the $\infty$\,-\,Ground State is $\infty$\,-\,harmonic. This explains how the two problems are connected! (The corresponding statement is false for a finite $p$.)  Indeed, from [LL3] we have:

    \begin{thm}\label{LLV}  Let $\boldsymbol\beta$ be a streamline of $v_{\infty}$ that is not an attracting one. Then it cannot meet any other streamline before it either meets (and joins) an attracting streamline or reaches the High Ridge. The speed $|\nabla v_{\infty}(\boldsymbol\beta(t))|$ is constant along $\boldsymbol\beta$ until it joins an attracting streamline. 
      
      In the open set $\Omega \setminus \Gamma $ the $\infty$\,-\,Ground State satisfies the $\infty$-Laplace  Equation $\Delta_{\infty} v_{\infty}\,=\,0$.
    \end{thm}

    We mention that $\mathrm{area}(\Gamma)\,=\,0$.
    
\section{Gradient convergence up to the boundary}
For our purposes, we need the convergence of the modulus of the gradient $\nabla v_p$  \emph{up to the boundary} of $\Omega$, as $p\to\infty$. A similar result is needed for the $p$-Potential.  We split the proof  for the two problems \eqref{ring} and \eqref{Lambda}.

 \subsection{The $p$\,-\,Potential Function} \label{sec:ppot}
We shall use the $p$-harmonic approximation $u_p \to u_{\infty}$ where
    \begin{equation}\label{pharmappr}
 \begin{cases}
    \Delta_{p}u_p\,&=\,0\quad \text{in} \quad G\,=\,\Omega \setminus K\\
   \quad u_p\,&=\,0 \quad \text{on} \quad \partial \Omega\\
  \quad  u_p\,&=\,1 \quad \text{on} \quad \partial K.\\
  \end{cases}
    \end{equation}
    As usual, $\Delta_pu\,=\,\nabla\!\cdot\!(|\nabla u|^{p-2}\nabla u)$. For $p > 2$ (in two dimensions) it is known that $u_p \in C(\overline{G})$ and that it takes the correct boundary values at each point. (This valuable property holds in arbitrary domains, be they convex or not.) We recall the following results of J. Lewis in [L]; see also [Ja]:
    \begin{itemize}
    \item $u_p \nearrow u_{\infty}$ in $\overline{G}$.
    \item $\nabla u_p \neq 0$ in $G$.
    \item $u_p$ is real-analytic in $G$.
    \item $u_p$ has convex level curves.
    \item $u_p$ is superharmonic in $G$.
    \end{itemize}
    We shall need continuous second  derivatives on the boundary $\partial \Omega$. If $\partial \Omega$ is of class $C^2$, then $|\nabla u_p|\geq\nu_p>0$ in $\Omega \setminus K$  according to Lemma 2 in [L].  It is known that $ u_p \in C^1(\overline{\Omega}\setminus K)$. By classical theory for the  equation
    $$\sum_{i,j}\Bigl(\delta_{i\,j}|\nabla u_p|^2 + (p-2)\frac{\partial u_p}{\partial x_i}\,\frac{\partial u_p}{\partial x_j}\Bigr)\frac{\partial ^2 w}{\partial x_i \partial x_j} \,=\,0$$
    with ``frozen coefficients''
    we can conclude, using the Calderon-Zygmund theory, that the particular solution $w=u_p \in C^{2,\alpha}(\overline{\Omega}\setminus K)$ provided that
    $\partial \Omega$ is of class $C^{2,\alpha}$.  See Theorem 6.14 in [GT]. This is sufficient for our purpose.
        
   When $p = \infty$  we have that $u_\infty \in C^1(\overline{\Omega} \setminus K)$ if $\partial\Omega$ is of class $C^2$, according to Theorem 1.1 in [WY]. By results in [KZZ], the convergence $|\nabla u_p| \to |\nabla u_{\infty}|$ holds  \emph{locally} uniformly in $G$. Notice the absolute values! See Section 5 and Theorem 7 in [LL2] for a clarification. We shall need the convergence also at the outer boundary.

    \begin{lemma}\label{graduconv} Let $\partial \Omega$ be of class $C^2$. Then 
$$
 \lim_{p\to \infty}|\nabla u_p(\xi)|\,=\,|\nabla u _{\infty}(\xi)|,\qquad \xi\in\partial\Omega.
   $$
\end{lemma}
      \bigskip

      \emph{Proof:}  Let $\mathbf{n}$ denote the outer unit normal at $\xi$. By Theorem 1.1 in  [WY]  $\nabla u_{\infty} \in C(\overline{\Omega} \setminus K)$. We see that
      $$-\,\frac{\partial u_p(\xi)}{\partial \mathbf{n}}\,=\,|\nabla u_p(\xi)|$$
      for $2 < p \leq \infty$. Since $u_p \nearrow u_{\infty}$ and $u_p(\xi) = u_{\infty}(\xi)$, we can for $h>0$ write
    $$  \frac{u_{\infty}(\xi - h\mathbf{n})-u_{\infty}(\xi)}{h}\,\geq \, 
      \frac{u_{p}(\xi - h\mathbf{n})-u_{p}(\xi)}{h}$$
      and let $h \to 0$ to obtain
      $$-\,\frac{\partial u_{\infty}(\xi)}{\partial \mathbf{n}}\,\geq -\,\frac{\partial u_p(\xi)}{\partial \mathbf{n}},\qquad -\,\frac{\partial u_{\infty}(\xi)}{\partial \mathbf{n}}\,\geq \lim_{p \to \infty}\Bigl(-\,\frac{\partial u_p(\xi)}{\partial \mathbf{n}}\Bigr).$$
      The limit exists by monotonicity.
      
      For the reverse inequality, choose an interior disk $|(\xi - \varepsilon \mathbf{n})-x| < \varepsilon$ tangent to $\partial \Omega$ at the point $\xi$. By the comparison principle in the ring
      $$0 < |(\xi - \varepsilon \mathbf{n})-x| < \varepsilon$$
      we have
      $$u_p(x)\,\geq\,\frac{\varepsilon^{(p-2)/(p-1)} - |x -(\xi - \varepsilon \mathbf{n})|^{(p-2)/(p-1)}}{\varepsilon^{(p-2)/(p-1)}}\,u_p(\xi - \varepsilon \mathbf{n}),$$
      where the fraction in the minorant is the fundamental solution of the $p$-Laplace Equation. At the point $\xi$ this inequality can be differentiated in the normal direction. It follows that
      $$   - \, \frac{\partial u_p(\xi)}{\partial \mathbf{n}}\, \geq \, \frac{p-2}{p-1}\,\frac{u_p( \xi - \varepsilon \mathbf{n})}{\varepsilon}$$
      and hence
      \[
\lim_{p\to \infty}\Bigl(-\,\frac{\partial u_p(\xi)}{\partial \mathbf{n}} \Bigr)\,\geq\, \frac{u_{\infty}( \xi - \varepsilon \mathbf{n})}{\varepsilon}
        =\, \frac{u_{\infty}( \xi - \varepsilon \mathbf{n})- u_{\infty}(\xi)}{\varepsilon}\,\to\, -\, \frac{\partial u_{\infty}(\xi)}{\partial \mathbf{n}},
      \]
	as $\varepsilon \to 0.$      This concludes the proof.\qquad $\Box$

 \subsection{The $p$-eigenvalue problem}\label{sec:peig}
The boundary convergence for $|\nabla v_p|$, where $v_p$ is the $p$-eigenfunction in equation (\ref{450}) has a  slightly different proof.

We shall later need continuous second derivatives on the boundary. Again, in a boundary zone, say $0<\mathrm{dist}(x,\partial \Omega) < \delta,$ there holds $|\nabla v_p|\geq \tilde \nu_p >0$, when $p$ is large. Indeed, since  $2 v_p\geq u_p$ on $K$ for large $p$ by \eqref{eq:normal}, where $u_p$ is the $p$-Potential with $K$ chosen to be the High Ridge (see \eqref{High}), the comparison principle implies that $2|\nabla v_p|\geq|\nabla u_p| \geq \nu_p$ at the boundary. By continuity it follows that near the boundary $|\nabla v_p|\geq \tilde \nu_p >0$ for some $\tilde\nu_p$.  Therefore, we may as in Section \ref{sec:ppot}  conclude that $v_p$ has continuous second order derivatives on the boundary, provided that $\partial \Omega$ is of class $C^{2,\alpha}$.

Near the boundary, $v_{\infty}$ is a solution to $\Delta_{\infty}v_{\infty}\,=\,0$, see Corollary 7 in [LL3]. In particular, $\nabla v_{\infty}$ is continuous up to the boundary in this zone, see Theorem 1.1 in [WY].

     \begin{lemma}\label{gradvconv} Let $\partial \Omega$ be of class $C^2$. Then 
     $$
 \lim_{p\to \infty}|\nabla v_p(\xi)|\,=\,|\nabla v _{\infty}(\xi)|,\qquad \xi \in \partial \Omega.
    $$
\end{lemma}
      \bigskip

      \emph{Proof:}  The inequality
      \begin{equation}\label{liminf}\liminf_{p \to \infty}|\nabla v_p(\xi)|\,\geq\, |\nabla v_p(\xi)|\end{equation}
      comes from a similar comparison as in Lemma \ref{graduconv}.  Now
      $$ \Delta_pv_p\,=\, - \lambda_pv_p^{p-1} \,<\,0 $$
      so that $v_p$ is a supersolution of the $p$-Laplace Equation. Let $\mathbf{n}$ denote the outer unit normal at $\xi$. Using the comparison principle in the interior ring $0 < |x -(\xi-(\varepsilon \mathbf{n})| < \varepsilon$ we have
      $$v_p(x)\,\geq\,\frac{\varepsilon^{(p-2)/(p-1)} - |x -(\xi - \varepsilon \mathbf{n})|^{(p-2)/(p-1)}}{\varepsilon^{(p-2)/(p-1)}}\,v_p(\xi - \varepsilon \mathbf{n}),$$
      where the fraction in the minorant is the fundamental solution of the $p$-Laplace Equation. Now the inequality (\ref{liminf}) comes as in the previous subsection.

      The reverse inequality requires some tinkering, because we do not know whether $v_{\infty} \geq v_p$. Consider the ascending streamline $\boldsymbol\beta_p = \boldsymbol\beta_p(t)$ for the $p$-eigenfunction $v_p$ starting at $\xi$:
      $$\frac{d\,\,}{dt}\, \boldsymbol\beta_p(t)\,=\, +\,\nabla v_p(\boldsymbol\beta_p(t)),\qquad 
      \boldsymbol\beta_p(0)\,=\,\xi .$$ (The \emph{ascending} streamlines are unique.)
        According to the end of Section 4 in [LL3] we have
        $$\frac{d\,\,}{dt}\Bigl(|\nabla v_p(\boldsymbol\beta_p(t))|^2 + \frac{1}{2}\kappa_p t^2\Bigr)  \,\geq\,0,\qquad t > 0,$$
        where the constant $\kappa_p \to 0+$ as $p \to \infty$. It follows that
        $$ |\nabla v_p(\boldsymbol\beta_p(t))|^2 + \frac{1}{2}\kappa_p t^2\,\geq\, |\nabla v_p(\boldsymbol\beta_p(0))|^2\,=\,  |\nabla v_p(\xi)|^2.$$
        By the results in [LL3] (see the proof of Theorem 10) we know that $\nabla v_p \to \nabla v_{\infty}$ \emph{locally} uniformly and that $\boldsymbol\beta_p\,\to\, \boldsymbol\beta_{\infty}$
        pointwise\footnote{In particular, the bound for $\|\nabla v_p\|_{\infty}$ in Lemma 5 in the arXiv version of [LL3] yields $$|\boldsymbol\beta_p(t_2)-\boldsymbol\beta_p(0)| = \left\vert \int_0^{t_2}\!\nabla v_p(\boldsymbol\beta_p(t))\,dt \right\vert \,\leq(\lambda_p\mathrm{diam}\Omega)^{\frac{1}{p-1}}\|v_p\|_{\infty}\cdot (t_2-0)$$
        and it follows that          $$|\boldsymbol\beta_{\infty}(t)-\xi|\,\leq\,\Lambda\|v_{\infty}\|_{\infty}\cdot t.$$}, where $\boldsymbol\beta_{\infty}$ is the streamline of $v_{\infty}$ emerging at $\xi$. We conclude that
        $$|\nabla v_{\infty}(\boldsymbol\beta_{\infty}(t))|\,\geq\, \limsup_{p \to \infty}|\nabla v_p(\xi)|,\qquad t>0.$$
        Since $\nabla v_{\infty}$ is continuous up to the boundary, by sending $t$ to $0+$ we finally arrive at
        $$|\nabla v_{\infty}(\xi) |\,\geq\, \limsup_{p \to \infty}|\nabla v_p(\xi)|.$$
        Thus the lemma is proved. \qquad  $\Box$

        \section{Assumptions for the Moving Plane Method}

        In the plane, Alexandroff's method is about a moving \emph{line}, across which solutions are reflected. For simplicity, we immediately make the following assumptions.

        \bigskip
        \textbf{Assumptions:} Suppose from now on that
          \begin{description}
          \item{1.}\quad $\Omega$ is a bounded convex domain in the $xy$-plane.
          \item{2.}\quad $\Omega$ is symmetric with respect to the $x$-axis and $y$-axis.
          \item{3.}\quad $\partial \Omega$ is of class $C^{3,\alpha}$.
          \item{4.}\quad The curvature of $\partial \Omega$ is non-decreasing in the first quadrant when $x$ increases.
          \item{5.}\quad $K$ is the origin.
          \end{description}

\bigskip

        See [Ka], [CFP], and [S2].  We note that the above assumptions are valid for the  case when $\Omega$ is an ellipse in proper position.
        
        We restrict our description to the first quadrant and consider a non-horizontal line $\ell$. The line  divides the plane in two open half-planes $T_+$ and $T_-$, where $T_-$ is chosen so that $T_-$ lies to the right of $\ell$. Let $x^*$ denote the reflexion of the point $x\in T_{-}$ across the line $\ell$. The above assumptions are designed to guarantee that  
        $$ (\Omega \cap T_-)^*\,\subset\,\Omega,$$
when $\ell$ is a normal to $\partial \Omega$. In other words, reflexion in the normals is possible.  See Figure \ref{reflfig}.  
        \begin{lemma}\label{possible} If the above assumptions are valid, then reflexion in the normals  is possible. In the first quadrant, the orientation is chosen so that $\Omega \cap T_-$ is to the right of the normal through a boundary point.
        \end{lemma}

        \emph{Proof:} See Lemma 4.2 in [CFP].\qquad    $\Box$

        \bigskip

        We define the reflected function of $f\in C(\overline{\Omega})$ as
        $$f^*:\,(\overline{\Omega} \cap T_-)^*\,\to\,\mathbb{R},\qquad f^*(x^*) = f(x).$$

   \paragraph{The $p$-eigenvalue problem.}  We use the $p$-eigenfunctions $v_p$ in equation (\ref{450}) and recall that they are continuous up to the boundary.
        \begin{lemma}\label{star} The reflected function $v_p^*$ satisfies
          $v_p\,\geq\,v_p^*$ in $(\Omega \cap T_-)^*.$
        \end{lemma}

        \bigskip
        \emph{Proof:}  Obviously  $v_p^*$ satisfies the same equation (\ref{450}) as $v_p$. Now $v_p\geq v_p^*$ on the boundary $\partial(\Omega \cap T_-)^*$. Indeed, $v_p = v_p^*$ on the line $\ell$ and $v_p \geq 0 =v_p^*$ on $(\partial \Omega)^*$. By the comparison principle $v_p \geq v_p^*$ in $(\Omega \cap T_-)^*$. The proof in [Li]  works for the comparison principle. \qquad  $\Box$
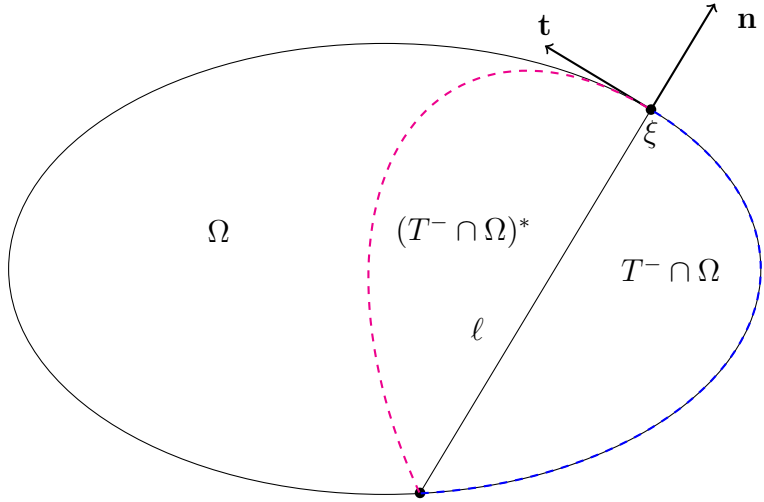
\begin{figure}
\def\originpath{
    plot[domain=-84.5:45, samples=60] ({5*cos(\x)-0.47},{3*sin(\x)+2.98}) 
}

\begin{center}
\begin{tikzpicture}[dot/.style={draw,fill,circle,inner sep=1pt}]
  \def\a{5} 
  \def\b{3} 
 
  \draw (0,0) ellipse ({\a} and {\b});
  \fill (3.54,2.12) circle[radius=2pt] node[below] {$\xi$};
    \fill (1,-0.8)  node[right] {$\ell$};
     \fill (-2.5,0.5)  node[right] {$\Omega$};
      \fill (0,0.5)  node[right] {$(T^-\cap \Omega)^*$};
       \fill (3,0)  node[right] {$T^-\cap \Omega$};
       \draw[->,thick] (3.54,2.12) -- (2.13,2.97);
       \fill (2.13,2.97) node[above] {$\bf t$};   
       \draw[->,thick,shift={(3.54,2.12)},rotate=-90] (0,0) -- (2.13-3.54,2.97-2.12);
       \fill (4.54,3.3) node[right] {$\bf n$};   
  \draw (3.54,2.12) -- (0.47,-2.98);
   \fill (0.47,-2.98) circle[radius=2pt];

\draw[blue,thick,dashed,shift={(0.47,-2.98)}] \originpath ;
\draw[magenta,thick,dashed,shift={(0.47,-2.98)},rotate=59,yscale=-1,rotate=-59] 
\originpath ;  
\end{tikzpicture}
\end{center}
\caption{The reflection illustrated in the case when $\Omega$ is an ellipse.}
\label{reflfig}
\end{figure}
        \begin{prop}\label{vtang} Let $\mathbf{n}$ denote the outer unit normal at the boundary point $\xi \in \partial \Omega$. Let $\mathbf{t}$ be the unit vector orthogonal to $\mathbf{n}$ and pointing from $T_-$ to $T_+$. Then we have
          at  all points  $\eta = \xi-\varepsilon \mathbf{n}$
          lying on the normal line in $\Omega$  that
          $$\frac{\partial v_p(\eta)}{\partial \mathbf{t}}\,\geq \,0.$$
        \end{prop}

        \emph{Proof:} Take $h>0$. By the previous lemma
        $$\frac{v_p(\eta+h\mathbf{t})-v_p(\eta)}{h}\,\geq\,\frac{v_p^*(\eta+h\mathbf{t})-v_p^*(\eta)}{h}\,=\,\frac{v_p(\eta-h\mathbf{t})-v_p(\eta)}{h}$$
        and as $h\,\to\,0$ we see that
        $$\frac{\partial v_p(\eta)}{\partial \mathbf{t}}\,\geq\, -\,\frac{\partial v_p(\eta)}{\partial \mathbf{t}}.$$
        This proves the desired inequality.\qquad $\Box$
           
        \bigskip

        For the next theorem we recall  that $v_p$ is of class $C^1(\overline{\Omega})$ and that $v_p$ is of class $C^2$ in a boundary zone $\overline{\Omega} \cap \{x|\,0<\mathrm{dist}(x,\partial \Omega) < \delta\}$. See Section \ref{sec:peig}. At  boundary points $\xi \in \partial \Omega$ we naturally have
        $$\frac{\partial v_p(\xi)}{\partial \mathbf{t}}\,=\,0 \quad\text{and}\quad \frac{\partial v_p(\xi)}{\partial \mathbf{n}}\,=\,-\,|\nabla v_p(\xi)|.$$

        \begin{thm}\label{thmvp} At the boundary point $\xi$ 
          $$\frac{\partial\,\,}{\partial \mathbf{t}}\,|\nabla v_p(\xi)|\,\geq\,0,$$
          where the tangent $\mathbf{t}$ at $\xi$ points from $T_-$ to $T_+.$
        \end{thm}

        \emph{Proof:} We first claim that
        $$\frac{\partial\,\,}{\partial \mathbf{n}}\Bigl(\frac{\partial v_p}{\partial \mathbf{t}}\Bigr)\,\leq\,0$$
        at $\xi$. To see this, let $\varepsilon\,\to\,0+$ in the difference quotient
        $$\frac{\dfrac{\partial v_p(\xi - \varepsilon \mathbf{n})}{\partial \mathbf{t}}\,-\,
          \dfrac{\partial v_p(\xi)}{\partial \mathbf{t}}}{\varepsilon}\,\geq\,0, $$
        which follows from Proposition \ref{vtang}, since $\tfrac{\partial v_p}{\partial \mathbf{t}} = 0$ at the boundary. To conclude the proof, use
        $$\frac{\partial\,\,}{\partial \mathbf{t}}\,|\nabla v_p|\,=\,-\,\frac{\partial\,\,}{\partial \mathbf{t}}\Bigl( \frac{\partial v_p}{\partial \mathbf{n}}\Bigr)\,=\,-\,\frac{\partial\,\,}{\partial \mathbf{n}}\Bigl( \frac{\partial v_p}{\partial \mathbf{t}}\Bigr).$$
        The mixed partial derivatives do not commute in general, but here they do, again due to $\tfrac{\partial v_p}{\partial \mathbf{t}} = 0$, so that the last term disappears from the general formula
        $$\frac{\partial\,\,}{\partial \mathbf{t}}\Bigl( \frac{\partial v_p}{\partial \mathbf{n}}\Bigr)\,=\frac{\partial\,\,}{\partial \mathbf{n}}\Bigl( \frac{\partial v_p}{\partial \mathbf{t}}\Bigr) + \kappa(\xi)\,
        \frac{\partial v_p}{\partial \mathbf{t}}  $$
        in differential geometry. See equation (10) in [Ka] or page 45 in [Sp].\qquad $\Box$
        
            \paragraph{The $p$\,-\,Potential Function.} For the $p$\,--\,potential function $u_p$ solving problem (\ref{pharmappr}) we encounter an extra problem caused by the inner boundary $\partial K$, which might hinder reflexion in some normals. Even if $(\Omega \cup T_-)^* \subset \Omega$, it is difficult to control $K$: often it happens that $(K \cup T_-)^*$ contains points in $\Omega \setminus K$, when the line $\ell$ of reflexion is a normal to $\partial \Omega$.             
           In order to avoid a detailed geometric description, we have therefore chosen to assume that $K$ is the origin.

        \begin{thm}\label{thmup} Assume that $\Omega$ satisfies the assumptions in Section 4 and that $K= \{(0,0)\}$. Then the inequality
          $$\frac{\partial \,\,}{\partial \mathbf{t}}\,|\nabla u_p(\xi)|\,\geq\,0$$
          is valid when $\xi \in \partial \Omega$, where the tangent $ \mathbf{t}$ at $\xi$ points from $T_-$ to $T_+$. (That is, in the direction of decreasing curvature.)
        \end{thm}

        \emph{Proof:} We follow the same steps as in the proof for $v_p$. First, the counterparts to Lemmas \ref{possible} and \ref{star} follow as before, when one notices that the presence of $K$ does not spoil the comparison $u_p \geq u_p^*$. Indeed, if   $\xi_2 = f(\xi_1)$ is the equation of the boundary $\partial \Omega$ in the first quadrant $\xi_1>0,\, \xi_2>0$, then the normal through $(\xi_1,\xi_2)\in \partial \Omega$ intersects the $x$-axis at the point
        $$\bigl(f(\xi_1)f'(\xi_1)+\xi_1 ,\,0\bigr).$$
        By Lemma 4.3 in [CFP], the assumptions imply that $f(\xi_1)f'(\xi_1)+\xi_1 \,\geq\,0$. In other words, $K$ (the origin) is not reflected at all. Thus nothing hinders the comparison $u_p \geq u_p^*$. This yields Lemma \ref{star} for $u_p$.
        The counterpart to Proposition \ref{vtang} follows. So does Theorem \ref{thmvp} for $u_p$.\qquad $\Box$

        \section{The Passage to $u_{\infty}$ and $v_{\infty}$}

        Assume again that the assumptions on the domain in Section 5 are fulfilled and that $K$ is the origin. By Theorems \ref{thmvp} and \ref{thmup}
        $$\frac{d \,\,}{d \mathbf{t}}\,|\nabla u_p(\xi)|\,\geq\,0,\qquad   \frac{d \,\,}{d \mathbf{t}}\,|\nabla v_p(\xi)|\,\geq\,0$$
        where the tangent $ \mathbf{t}$ points in the direction of non-increasing curvature. That is to the left in the first quadrant. If now
        $\xi$ and $\zeta$ belong to $\partial \Omega$ in the first quadrant and are ordered  so that  $\xi$ is to the left of $\zeta$, then
        $$|\nabla u_p(\xi)|\,\geq\, |\nabla u_p(\zeta)|\quad\text{and}\quad 
          |\nabla v_p(\xi)|\,\geq\, |\nabla v_p(\zeta)|$$
        i.e.,  $ |\nabla u_p|$  and $|\nabla v_p|$ increase when the curvature decreases.

        On the boundary,  Lemma \ref{graduconv} and Lemma \ref{gradvconv} assure that we can proceed to the limits. We arrive at
        
       \begin{equation}\label{monotone}    
          |\nabla u_{\infty}(\xi)|\,\geq\, |\nabla u_{\infty}(\zeta)|\quad\text{and}\quad 
                |\nabla v_{\infty}(\xi)|\,\geq\, |\nabla v_{\infty}(\zeta)|.
       \end{equation}

       This monotonicity enables us to conclude that  $|\nabla u_{\infty}|$ has only two maxima and two minima on the  boundary $\partial \Omega$, viz. at the four intersections with the coordinate axes: the maxima are on the $y$-axis, the minima on the $x$-axis. ---The same goes for  $|\nabla v_{\infty}|$.

       We want to show that this monotonicity is strict. To see this, assume that $|\nabla u_{\infty}(\xi)| = |\nabla u_{\infty}(\zeta)|$ for two different boundary points in the first quadrant. Then $|\nabla u_{\infty}|$ would be constant along the boundary arc between the points. According to Lemma 12 and Lemma 16 in [LL2]  this means that all the streamlines emerging from this arc are \emph{straight line segments} that cannot intersect each others, except at $K$, which now is the origin.  But this forces the boundary arc to be an arc of a circle centered at the origin.
       
       The minimum of $|\nabla u_{\infty}|$ on the boundary is attained at the $x$-axis. We claim that the minimum is strict. If not, we would have $|\nabla u_{\infty}| = c$ on a circular arc. The above mentioned Lemma 12 and Lemma 16 in [LL2] also imply that the eikonal equation  $|\nabla u_{\infty}(x,y)| = c$ is valid in the whole closed circular sector. Assuming that $\Omega$ is not a disk, which case is trivial, we choose a boundary point $\xi$ not on the circular arc. Then $|\nabla u_{\infty}(\boldsymbol\alpha(0)| = C >c$ .
      
       Recall that the (ascending) streamlines $\boldsymbol\alpha = \boldsymbol\alpha(t)$ are defined through
       $$\frac{d \boldsymbol\alpha(t)}{d t}\,=\,+\,\nabla u_{\infty}(\boldsymbol\alpha(t)), \qquad \boldsymbol\alpha(0)\,=\,\xi.$$
       They start at the boundary and reach the origin. Always, the speed $|\nabla u_{\infty}(\boldsymbol\alpha(t)|$ is non-decreasing. We see that 
      $$\limsup_{(x,y)\to (0,0)}|\nabla u_{\infty}(x,y)|\,\geq\,C,\quad \liminf_{(x,y)\to (0,0)}|\nabla u_{\infty}(x,y)|\,\leq\,c.$$
       This contradicts Proposition 10  in [LL1] according to which the full limit exists at the origin. Therefore \emph{the minimum is strict.}

       The two streamlines starting at the intersection of $\partial \Omega$ with the $x$-axis are \emph{attracting streamlines} in the terminology of [LL2]. By symmetry, they are line segments on the $x$-axis.
        Now Theorem \ref{LLU} can be stated in the following form. 

       \begin{prop}\label{mainu}  Suppose that the assumptions in Section 4 are valid and assume that the domain is not a disk.  Let $\boldsymbol \alpha$ be a  streamline whose initial point is not on the $x$-axis.  It cannot meet any other streamline before it meets and joins the $x$-axis. The speed $|\nabla u_{\infty}(\boldsymbol \alpha(t)|$ is constant along $\boldsymbol \alpha$ until it meets   the $x$-axis,  after which the speed is non-decreasing.
        \end{prop}

       \emph{Proof:}  Equation (\ref{monotone}) and the above discussion allows us to conclude this from Theorem 3 in [LL2].  \qquad $\Box$

       \bigskip
       
       A similar version of Theorem \ref{LLV} holds for $v_{\infty}$. Now circular boundary arcs where $|\nabla v_{\infty}|$ is constant are not excluded. In addition, we can infer the following interesting property.

       \begin{prop} Suppose that the assumptions in Section 4 are valid. Then
         the $\infty$\,-\,Ground State satisfies the equation $$\Delta_{\infty}v_{\infty}\,=\,0\quad \text{in}\quad \Omega\quad\text{except possibly on the}\, x\!-\!\text{axis.}$$
         Streamlines cannot meet outside the $x$-axis.
       \end{prop}

       \emph{Proof:} This  essentially follows from Theorem \ref{LLV}. To see this, we first claim that the minimum boundary speed $|\nabla v_{\infty}|$, which is attained  at the $x$-axis, is either strict or  is attained along a circular arc. The strict case is immediately clear by Theorem \ref{LLV}.
       
       In order to treat the other case, we need  the contact set $\Upsilon$, see Section \ref{sec:prel}.
             It is closed, has zero area and $v_\infty$ is $C^1$ and satisfies $\Delta_{\infty} v_{\infty} =0$ outside $\Upsilon$. 
       Moreover, $\Upsilon$ does not touch $\partial \Omega$.
       
       Suppose now that we have a  closed boundary arc $C$ which is symmetric about the $x$-axis  and that the speed is constant along it. We can assume that it is of maximal length: the boundary speed outside is strictly larger.
       By symmetry, the two streamlines $\boldsymbol \gamma_1$ and $\boldsymbol \gamma_2$, starting at its  endpoints, intersect at a point $P$ on the $x$-axis. We shall now argue that $P$ is the only point of  the contact set  $\Upsilon$ lying in the \emph{closed} region  bounded by $\boldsymbol\gamma_1$, $\boldsymbol \gamma_2$ and the boundary arc $C$.  Suppose, towards a contradiction, that the lowest level curve $\boldsymbol \omega$ that in  this closed region  reaches $\Upsilon$ does not contain $P$. (Thus $P$ is at a higher level.) So  $v_\infty$ is of class $C^1$ and satisfies $\Delta_{\infty} v_{\infty} =0$ in the  open region
       bounded by $\boldsymbol\omega$, $\boldsymbol\gamma_1$, $\boldsymbol\gamma_2$ and $C$. Therefore, we can apply Lemma 12 in [LL2] to conclude that the
       $$\text{eikonal equation}\quad |\nabla v_{\infty}| = c$$
       is valid in this region.
       
       Along any lower level curve, say $\boldsymbol\omega_-$, we now have that
       $$\frac{|\nabla v_{\infty}(\boldsymbol\omega_-)|}{|v_{\infty}(\boldsymbol\omega_-)|}\,= \, \frac{ c}{v_{\infty}(\boldsymbol\omega_-)}\, =\quad \text{constant}.$$
       The approximation $\boldsymbol\omega_- \to \boldsymbol\omega$ of the level curve from below implies, by the continuity of $S^-$ operator, that the whole arc of the level curve $ \boldsymbol\omega$ between $\boldsymbol \gamma_1$ and $\boldsymbol \gamma_2$   belongs to $\Upsilon$. By Lemma 9 in [LL3], the whole sector between  $\boldsymbol\omega$, $\boldsymbol\gamma_1$ and $\boldsymbol\gamma_2$ (with apex at $P$) belongs to $\Upsilon$. This contradicts the fact that $\Upsilon$ has zero area. 
        
       Hence, the first point in $\Upsilon$ is $P$. In particular, the streamlines $\boldsymbol\gamma_1$ and $\boldsymbol\gamma_2$ do not contain any points of $\Upsilon$ below $P$.  Hence, $v_\infty$ is of class $C^1$ and satisfies $\Delta_{\infty} v_{\infty} =0$ in the whole  region bounded by $\boldsymbol\gamma_1$, $\boldsymbol\gamma_2$ and $C$.   Again by Lemma 12 in [LL2], the eikonal equation holds here. It follows from Lemma 1 in [Ar2]  that the streamlines emerging from the arc $C$ are non-intersecting straight lines intersecting only at the point $P$. This implies that the arc $C$ has to be circular.
       
       Moreover, since $|\nabla v_{\infty}|$ is constant along $\boldsymbol\gamma_1$ and $\boldsymbol\gamma_2$ until they meet at $P$, no streamline emerging from the part of $\partial \Omega$ that is outside of $C$ can meet $\boldsymbol\gamma_1$ or $\boldsymbol\gamma_2$ below $P$, because the emerging streamline has too high an initial speed.
                 
       In conclusion, if the minimum is not strict we have a \emph{circular} boundary arc with constant speed, and the streamlines are rays joining at a point on the $x$-axis.\qquad $\Box$
       \bigskip
       
       In fact, one can extract more, but we are content to provide one good example.

       \bigskip
      \noindent  \textbf{Example:}  The ellipse
       $$\frac{x^2}{a^2}+\frac{y^2}{b^2}\,=\,1,\qquad  0 < b < a,$$
       fulfills all our assumptions. Given the standard parametrization $x(t)=a\cos t$ and $y(t)=b\sin t$, the curvature at $(x(t),y(t))$ is given by
       $$
       \frac{ab}{\left(a^2\sin^2t+b^2\cos^2 t\right)^\frac32}.
       $$
       Clearly the curvature is decreasing from $t=0$ to $t=\pi/2$. The rest of the assumptions in Section 4 are obviously satisfied. ---See also page 215 in[CFP].

       \section{The $\infty$\,-\,Potential in a Square}\label{sec:brustad}

       The remarkable formula
       \begin{equation*}
u_{\infty}(x,y)\,=\,\min_{0\leq\psi\leq\frac{\pi}{2}}\max_{0\leq \rho\leq 1}\left\{x\,\rho\cos(\psi) + y\, \rho \sin(\psi)-W(\rho,\psi) \right\}
       \end{equation*}
       where the function $W(\rho,\psi)$ has the explicit representation
       $$ W(\rho,\psi)\,=\,\frac{8}{\pi}\Bigl(\frac{\rho^4}{6}\sin(2\psi)+\frac{\rho^{36}}{210}\sin(6\psi) +\frac{\rho^{100}}{990}\sin(10\psi) +\cdots \Bigr) \nonumber
         $$
       was discovered in [B] for the $\infty$-Potential $u_{\infty}$  of the punctured square $0<|x-1|<1,\,\,0 < |y-1|< 1$. Here the center $(1,1)$ is removed from the square. The formula is valid in the subsquare $0 \leq x \leq 1,\,\,0 \leq y \leq 1$ and is extended by symmetry. The resulting function is of class $C^1$ up to the sides (but not at the center) and it is real-analytic outside the diagonals $y-1 = \pm (x-1)$. See Figure \ref{fig:inftypot}.
\begin{center}
\begin{figure}[h!]
\begin{center}
\includegraphics[scale=1.3]{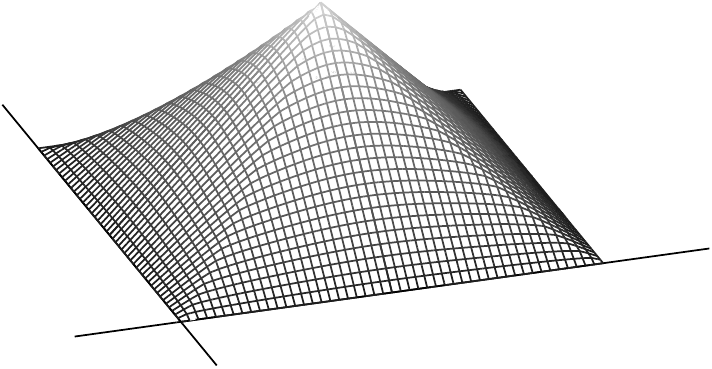}
\end{center}
\caption{The $\infty$-Potential in the square.}
\label{fig:inftypot}
\end{figure}
\end{center} 
       \emph{This function is not equal to the $\infty$-Ground State} $v_{\infty}$. See Section 5 in [B] for the original proof of this fact. Below, we briefly explain why they do not coincide.

       Observe that if the functions would coincide, then also $v_{\infty}$ should be of class $C^1$. By \eqref{Lambda}, 
       this implies        that
       $$\frac{|\nabla v_{\infty}|}{v_{\infty}}\,\geq \,1\quad\text{and hence}\quad 
       \frac{|\nabla u_{\infty}|}{u_{\infty}}\,\geq \,1$$
       in the punctured square. (Here $\Lambda =1$.) In particular, the last inequality should  hold on the diagonal. A numerical calculation below will show that this is not the case.

       From [B] we have, using the variables   
       \begin{equation}\label{eq:trho}
       t\,=\,\frac{8}{\pi}\,\sum_{n=1}^{\infty}(-1)^{n-1}\frac{m_n}{m_n^2-1}\,\rho^{m_n^2-1},\qquad m_n= 4n-2,
       \end{equation}
       that on the diagonal the so-called Rayleigh quotient takes the form
       $$R(t)\,=\,\frac{|\nabla u_{\infty}(t/\sqrt{2},t/\sqrt{2})|}{u_{\infty}(t/\sqrt{2},t/\sqrt{2})}\,=\,\frac{\rho}{\rho\, W_{\rho}(\rho,\frac{\pi}{4}) - W(\rho,\frac{\pi}{4})}.$$
       Here $0\leq t \leq \sqrt{2}$ and $0\leq \rho \leq 1$. The quantities involved are
       \begin{align*}
         W(\rho,\frac{\pi}{4})\,&=\,\frac{8}{\pi}\,\sum_{n=1}^{\infty}\frac{(-1)^{n-1}}{(m_n^2-1)m_n}\,\rho^{m_n^2}\\
         \rho\,W_{\rho}(\rho,\frac{\pi}{4})\,&=\,\frac{8}{\pi}\,\sum_{n=1}^{\infty}\frac{(-1)^{n-1}m_n}{m_n^2-1}\,\rho^{m_n^2}\\
         \rho\,W_{\rho}(\rho,\frac{\pi}{4})-W(\rho,\frac{\pi}{4})  \,&=\,\frac{8}{\pi}\,\sum_{n=1}^{\infty}\frac{(-1)^{n-1}}{m_n}\,\rho^{m_n^2}
       \end{align*}
       and so we arrive at the simple expression
       \begin{equation}\label{Req}
         \frac{1}{R}\,=\,\frac{8}{\pi}\,\sum_{n=1}^{\infty}\frac{(-1)^{n-1}}{m_n}\,\rho^{m_n^2-1
           }
       \end{equation} 
       for the reciprocal value of the  Rayleigh quotient. The series converges for $0\leq \rho \leq 1$, and its sum is $1$ when $\rho = 1$. To conclude, we only have to exhibit a value of $\rho$ for which $R >1$. When $\rho = 0.97$ we have by \eqref{Req}
       \[
       \begin{split}
       \frac{1}{R}\,&=\, \frac{8}{\pi}\left(\frac{1}{2} 0.97^3-\frac{1}{6}0.97^{35}  + \frac{1}{10}0.97^{99} - \cdots \right)\\
       &>\frac{8}{\pi}\left(\frac{1}{2} 0.97^3-\frac{1}{6}0.97^{35} \right)\\
       &\approx 1.01590>1.
       \end{split}
       \]
       Here we have used Leibniz's rule for alternating series. By \eqref{eq:trho}, $\rho=0.97$ corresponds to $t\approx 1.4112$. The corresponding point $(0.9979,\phantom{.\!}0.9979)$ on the diagonal is at the distance $0.0030$ from the center. 

       We can extract further information about the ``unknown'' variational\footnote{Recall that this is obtained as the limit of $v_{p}$. It is not known to be unique, but it inherits the symmetries of $v_p$.} $\infty$\,-\,Ground State $v_{\infty}$, be it unique or not, using the fact that it is not the $\infty$\,-\,Potential. First, it cannot be of class $C^1$ in the whole punctured square, because Theorem 3.1 in [Y] would then imply that the functions coincide. Second, using Theorem 2 we can deduce that the  variational $\infty$\,-\,Ground State $v_{\infty}$ is $\infty$-harmonic except on a portion of the diagonals lying in a symmetric neighbourhood around the center. In other words, the \emph{contact set} looks like the letter \textsf{X}, where the crossing line segments have length at most $2(\sqrt2-1)$; probably much shorter.
       \bigskip

       \section{Epilogue}
       The described results are based on a fairly recent theory. So several immediate questions seem to be open problems. A few of them are
       \begin{itemize}
       \item  For which convex domains do $|\nabla u_{\infty}(\xi)|$ or $|\nabla v_{\infty}(\xi)|$ have only a finite number of minima on the boundary?
         \item  Are $u_{\infty}$ and $v_{\infty}$ twice  differentiable or even real-analytic outside the attracting streamlines? 
       \item  Is $\log(v)$ concave for all solutions $v$ of the $\infty$\,-\,Eigenvalue Problem (\ref{Lambda}), be they variational or not?
       \item  Are there other domains than the stadiums in which we have  $u_{\infty}\,\equiv\,v_{\infty}$?
       \item  How do the streamlines run in non-convex domains?
         \item  What about several dimensions?
       \end{itemize}
       There are many more interesting questions, but we must stop here.
       \bigskip

{\small \noindent {\textsf{Karl Brustad\\Frostavegen 1691\\NO--7633 Frosta, Norway }}\\ 
  \textsf{e-mail}: brustadkarl@gmail.com\\

{\small \noindent {\textsf{Erik Lindgren\\  Department of Mathematics, KTH -- Royal Institute of Technology\\ 100 44, Stockholm, Sweden}}  \\
\textsf{e-mail}: eriklin@kth.se\\

{\small \noindent {\textsf{Peter Lindqvist\\ Department of
    Mathematical Sciences, Norwegian University of Science and
    Techno\-logy, NO--7491, Trondheim, Norway}}\\
\textsf{e-mail}: peter.lindqvist@ntnu.no

\end{document}